\documentclass[12pt]{article}
\usepackage[psamsfonts]{amssymb}
\usepackage{amsthm,amsmath}

\title{Axiomatizations of quasi-polynomial functions on bounded chains}

\author{Miguel Couceiro \\
Mathematics Research Unit, University of Luxembourg \\
162A, avenue de la Fa\"{\i}encerie, L-1511 Luxembourg, Luxembourg \\
miguel.couceiro[at]uni.lu
\and %
Jean-Luc Marichal \\
Mathematics Research Unit, University of Luxembourg \\
162A, avenue de la Fa\"{\i}encerie, L-1511 Luxembourg, Luxembourg \\
jean-luc.marichal[at]uni.lu }

\date{Revised, May 25, 2009}

\begin{document}
\maketitle

\theoremstyle{plain}
\newtheorem{theorem}{Theorem}
\newtheorem{lemma}[theorem]{Lemma}
\newtheorem{proposition}[theorem]{Proposition}
\newtheorem{corollary}[theorem]{Corollary}
\newtheorem{fact}[theorem]{Fact}
\newtheorem{claim}{Claim}

\theoremstyle{definition}
\newtheorem{definition}[theorem]{Definition}
\newtheorem{example}[theorem]{Example}

\theoremstyle{remark}
\newtheorem{conjecture}{Conjecture}
\newtheorem{remark}{Remark}

\newcommand{\card}[1]{\ensuremath{\lvert{#1}\rvert}}
\newcommand{\vect}[1]{\ensuremath{\mathbf{#1}}} 
\newcommand{\co}[1]{\ensuremath{\overline{#1}}}
\def\median{\mathop{\rm med}\nolimits}

\begin{abstract}
Two emergent properties in aggregation theory are investigated, namely horizontal maxitivity and comonotonic maxitivity (as well as their dual
counterparts) which are commonly defined by means of certain functional equations. We completely describe the function classes axiomatized by
each of these properties, up to weak versions of monotonicity in the cases of horizontal maxitivity and minitivity. While studying the classes
axiomatized by combinations of these properties, we introduce the concept of quasi-polynomial function which appears as a natural extension of
the well-established notion of polynomial function. We give further axiomatizations for this class both in terms of functional equations and
natural relaxations of homogeneity and median decomposability. As noteworthy particular cases, we investigate those subclasses of quasi-term
functions and quasi-weighted maximum and minimum functions, and provide characterizations accordingly.
\end{abstract}


\noindent{\bf MSC classes:} 28B15, 39B72 (Primary) 06A05, 06D05 (Secondary)

\section{Introduction}

Aggregation functions arise wherever aggregating information is important: applied and pure mathematics (probability, statistics, decision
theory, functional equations), operations research, computer science, and many applied fields (economics and finance, pattern recognition and
image processing, data fusion, etc.). For recent references, see Beliakov et al.~\cite{BelPraCal07} and Grabisch et al.~\cite{GraMarMesPap09}.

A noteworthy aggregation function is the so-called discrete Sugeno integral, which was introduced by Sugeno~\cite{Sug74,Sug77} and which has
been widely investigated in aggregation theory, due to its many applications for instance in fuzzy set theory, decision making, and image
analysis. For general background, see also the edited book \cite{GraMurSug00}.

A convenient way to introduce the discrete Sugeno integral is via the concept of (lattice) polynomial functions, i.e., functions which can be
expressed as combinations of variables and constants using the lattice operations $\wedge$ and $\vee$. As shown in \cite{Marc}, the
discrete Sugeno integrals are exactly those polynomial functions $f:L^n\to L$ which are idempotent, that is, satisfying $f(x,\ldots,x)=x$.
Several axiomatizations of the class of discrete Sugeno integrals (as idempotent polynomial functions) have been recently given; see
\cite{CouMar2}.

Of particular interest in aggregation theory are the so-called horizontal maxitivity and comonotonic maxitivity (as well as their dual
counterparts), usually expressed in terms of certain functional equations, and which we now informally describe.

Let $L$ be a bounded chain. For every $\vect{x}\in L^n$ and every $c\in L$, consider the horizontal maxitive decomposition of $\vect{x}$
obtained by ``cutting'' it with $c$, namely $\vect{x} = (\vect{x}\wedge c)\vee [\vect{x}]_c$, where $[\vect{x}]_c$ is the $n$-tuple whose $i$th
component is $0$, if $x_i\leqslant c$, and $x_i$, otherwise. A function $f\colon L^{n}\rightarrow L$ is said to be \emph{horizontally maxitive}
if
$$
f(\vect{x}) = f(\vect{x}\wedge c)\vee f([\vect{x}]_c)
$$
for every $\vect{x}\in L^n$ and every $c\in L$.

A function $f\colon L^{n}\rightarrow L$ is said to be
 \emph{comonotonic maxitive} if, for any two vectors $\vect{x}$ and $\vect{x'}$ in the same standard simplex of $L^n$, we have
$$
f(\vect{x}\vee \vect{x'}) = f(\vect{x})\vee f(\vect{x'}).
$$

As we are going to see (Lemma~\ref{lemma:ComHor} below), these properties (as well as their duals) are closely related and constitute properties
shared by discrete Sugeno integrals. Still, and as it will become evident, no combination of these with their dual forms suffices to fully
describe the class of Sugeno integrals. Thus, and given their emergence in aggregation theory, it is natural to ask which classes of functions
are axiomatized by combinations of these properties or, in fact, by each of these properties.

In this paper, we answer this question for both the maxitive and minitive comonotonic properties, and for horizontal maxitivity and minitivity
properties, up to certain weak variants of monotonicity. While looking at combinations of the latter properties, we reach a natural
generalization of polynomial functions, which we call \emph{quasi-polynomial functions} and which are best described by the following equation
\begin{equation*}
f(x_1,\ldots,x_n)=p(\varphi(x_1),\ldots,\varphi(x_n)),
\end{equation*}
where $p$ is a polynomial function and $\varphi$ a nondecreasing function (see Theorem~\ref{thm:QuasiPol} below). Searching for alternative
descriptions, we introduce weaker versions of well-established properties, such as homogeneity and median decomposability, to provide further
axiomatizations of the class of quasi-polynomial functions, accordingly.

This paper is organized as follows. We start by recalling basic notions and terminology in lattice  function theory, as well as present some
known results, needed throughout this paper (Section 2). In Section 3, we study the properties of horizontal maxitivity and comonotonic
maxitivity, as well as their dual forms, and determine those function classes axiomatized by each of these properties. Combinations of the
latter are then considered in Section 4.1, where the notion of quasi-polynomial function is introduced. In Section 4.2, we propose weaker
versions of homogeneity and median decomposability, and provide further characterizations of quasi-polynomial functions, accordingly. In Section
5, we introduce and axiomatize few noteworthy subclasses of quasi-polynomial functions, namely, those of quasi-term functions and those
quasi-weighted maximum and minimum functions.

\section{Basic notions and preliminary results}

In this section we recall basic terminology as well as some results needed in the current paper. For general background we refer the reader to,
e.g., Burris and Sankappanavar~\cite{BurSan81} and Rudeanu~\cite{Rud01}.

\subsection{General background}

Throughout this paper, let $L$ be a bounded chain with operations $\wedge$ and $\vee$, and with least and greatest elements $0$ and $1$,
respectively. A subset $S$ of a chain $L$ is said to be \emph{convex} if for every $a,b\in S$ and every $c\in L$ such that $a\leqslant
c\leqslant b$, we have $c\in S$. For any subset $S\subseteq L$, we denote by $\co{S}$ the convex hull of $S$, that is, the smallest convex
subset of $L$ containing $S$. For every $a,b\in S$ such that $a\leqslant b$, the \emph{interval} $[a,b]$ is the set $[a,b]=\{c\in L: a\leqslant
c\leqslant b\}.$ For any integer $n\geqslant 1$, let $[n]=\{1,\ldots,n\}$.

For any bounded chain $L$, we regard the Cartesian product $L^n$, $n\geqslant 1$, as a distributive lattice endowed with the operations $\wedge$
and $\vee$ given by
\begin{eqnarray*}
(a_1,\ldots ,a_n)\wedge (b_1,\ldots ,b_n) &=& (a_1\wedge b_1, \ldots , a_n\wedge b_n), \\
(a_1,\ldots ,a_n)\vee (b_1,\ldots ,b_n) &=& (a_1\vee b_1, \ldots , a_n\vee b_n).
\end{eqnarray*}
The elements of $L$ are denoted by lower case letters $a,b,c,\ldots$, and the elements of $L^n$, $n>1$, by bold face letters
$\vect{a},\vect{b},\vect{c},\ldots$. We also use $\vect{0}$ and $\vect{1}$ to denote the least element and greatest element, respectively, of
$L^n$. For $c\in L$ and $\vect{x}=(x_1,\ldots ,x_n)\in L^n$, set
$$
\vect{x}\wedge c=(x_1\wedge c, \ldots , x_n\wedge c) \quad \textrm{and} \quad \vect{x}\vee c=(x_1\vee c, \ldots , x_n\vee c).
$$

The \emph{range} of a function $f\colon L^{n}\rightarrow L$ is defined by $\mathcal{R}_f=\{f(\vect{x}) : \vect{x}\in L^n\}$. A function $f\colon
L^{n}\rightarrow L$ is said to be \emph{nondecreasing} (\emph{in each variable}) if, for every $\vect{a}, \vect{b}\in L^n$ such that
$\vect{a}\leqslant\vect{b}$, we have $f(\vect{a})\leqslant f(\vect{b})$. The \emph{diagonal section} of $f$, denoted $\delta_f$, is defined as
the unary function given by $\delta_f(x)=f(x,\ldots ,x)$. Note that if $f$ is nondecreasing, then $\delta_f$ is nondecreasing and
$\co{\mathcal{R}}_{\delta_f}=\co{\mathcal{R}}_f=[f(\vect{0}),f(\vect{1})]$.

\subsection{Polynomial functions and their representations}

In this paper the so-called polynomial functions will play a fundamental role. Formally, an \emph{$n$-ary polynomial function} on $L$ is any
function $f:L^n\to L$ which can be obtained by finitely many applications of the following rules:
\begin{enumerate}
\item[(i)] For each $i\in [n]$ and each $c\in L$, the projection $\vect{x}\mapsto x_i$ and the constant function $\vect{x}\mapsto c$ are
polynomial functions from $L^n$ to $L$.

\item[(ii)] If $f$ and $g$ are polynomial functions from $L^n$ to $L$, then $f\vee g$ and $f\wedge g$ are polynomial functions from $L^n$ to
$L$.
\end{enumerate}
Polynomial functions are also called lattice functions (Goodstein~\cite{Goo67}), algebraic functions (Burris and Sankappanavar~\cite{BurSan81})
or weighted lattice polynomial functions (Marichal~\cite{Marc}). Idempotent polynomial functions (i.e., satisfying $f(c,\ldots ,c)=c$ for every
$c\in L$) are referred to by aggregation theorists as (\emph{discrete}) \emph{Sugeno integrals}, and those obtained from projections by finitely
many applications of (ii) are usually referred to as (\emph{lattice}) \emph{term functions}. As an example, we have the ternary \emph{median function} \begin{eqnarray*}
\median(x,y,z) &=& (x\vee y)\wedge (y\vee z)\wedge (z\vee x)\\
&=& (x\wedge y)\vee (y\wedge z)\vee (z\wedge x).
\end{eqnarray*}

As observed by Goodstein~\cite{Goo67} (see also Rudeanu~\cite{Rud01}),  polynomial functions are exactly those functions which can be
represented by formulas in disjunctive and conjunctive normal forms. In fact, each polynomial function $f\colon L^{n}\rightarrow L$ is uniquely
determined by its restriction to $\{0,1\}^n$. Due to their relevance in the sequel, we recall some known results concerning normal form
representations of polynomial functions in the special case where $L$ is a chain. The following result is due to Goodstein~\cite{Goo67}.

\begin{proposition}\label{Goodstein}
\begin{itemize}
\item[(a)] Every polynomial function is completely determined by its restriction to $\{0,1\}^n$.

\item[(b)] A function $g\colon \{0,1\}^n\rightarrow L$ can be extended to a polynomial function $f\colon L^{n}\rightarrow L$ if and only if it
is nondecreasing. In this case, the extension is unique.

\item[(c)] For any $f\colon L^{n}\rightarrow L$, the following are equivalent:
\begin{enumerate}
\item[(i)] $f$ is a polynomial function.

\item[(ii)] There exists $\alpha \colon 2^{[n]}\rightarrow L$ such that
\begin{equation}\label{eq:DNF45}
f(\vect{x})=\bigvee_{I\subseteq [n]}\big(\alpha(I)\wedge \bigwedge_{i\in I} x_i\big).
\end{equation}
\item[(iii)] There exists $\beta \colon 2^{[n]}\rightarrow L$ such that
\begin{equation}\label{eq:CNF45}
f(\vect{x})=\bigwedge_{I\subseteq [n]}\big(\beta(I)\vee \bigvee_{i\in I} x_i\big).
\end{equation}
\end{enumerate}
\end{itemize}
\end{proposition}

The expressions given in (\ref{eq:DNF45}) and (\ref{eq:CNF45}) are usually referred to as the \emph{disjunctive normal form} (DNF)
representation and the \emph{conjunctive normal form} (CNF) representation, respectively, of the polynomial function $f$.

\begin{remark}
By requiring $\alpha$ and $\beta$ to be nonconstant functions from $2^{[n]}$ to $\{0,1\}$ and satisfying $\alpha(\varnothing)=0$ and
$\beta(\varnothing)=1$, respectively, we obtain the analogue of $(c)$ of Proposition~\ref{Goodstein} for term functions.
\end{remark}

As observed in \cite{Marc}, the DNF and CNF representations of polynomial functions $f\colon L^{n}\rightarrow L$ are not necessarily
unique. For instance, we have
$$
x_1\vee (x_1\wedge x_2)=x_1=x_1\wedge (x_1\vee x_2).
$$
However, from among all the possible set functions $\alpha$ (resp.\ $\beta$) defining the DNF (resp.\ CNF) representation of $f$, only one is
isotone (resp.\ antitone), namely the function $\alpha_f\colon 2^{[n]}\to L$ (resp.\ $\beta_f\colon 2^{[n]}\to L$) defined by
\begin{equation}\label{eq:AlphaBeta}
\alpha_f(I)=f(\vect{e}_I)\qquad\mbox{(resp.\ $\beta_f(I)=f(\vect{e}_{[n]\setminus I})$)},
\end{equation}
where $\vect{e}_I$ denotes the element of $\{0,1\}^n$ whose $i$th component is $1$ if and only if $i\in I$.

In the case when $L$ is a chain, it was shown in  \cite{CouMar2} that the DNF and CNF representations of polynomial functions $f\colon
L^{n}\rightarrow L$ can be refined and given in terms of standard simplices of $L^n$. Let $\sigma$ be a permutation on $[n]$. The \emph{standard
simplex} of $L^n$ associated with $\sigma$ is the subset $L^n_\sigma\subset L^n$ defined by
$$
L^n_\sigma =\{(x_1,\ldots ,x_n)\in L^n\colon x_{\sigma (1)}\leqslant x_{\sigma (2)}\leqslant \cdots\leqslant x_{\sigma (n)}\}.
$$
For each $i\in[n]$, define $S^{\uparrow}_\sigma (i)=\{\sigma(i),\ldots ,\sigma (n)\}$ and $S^{\downarrow}_\sigma (i)=\{\sigma(1),\ldots ,\sigma
(i)\}$. As a matter of convenience, set $S^{\uparrow}_\sigma (n+1)=S^{\downarrow}_\sigma (0)=\varnothing$.

\begin{proposition}\label{SimplexDNF}
For any function $f\colon L^{n}\rightarrow L$, the following conditions are equivalent:
\begin{enumerate}
\item[(i)] $f$ is a polynomial function.

\item[(ii)] For any permutation $\sigma$ on $[n]$ and every $\vect{x}\in L^n_\sigma$, we have
$$
f(\vect{x}) = \bigvee_{i=1}^{n+1}\big (\alpha_f(S^{\uparrow}_\sigma (i))\wedge x_{\sigma(i)}\big ),
$$
where $x_{\sigma(n+1)}=1$.

\item[(iii)] For any permutation $\sigma$ on $[n]$ and every $\vect{x}\in L^n_\sigma$, we have
$$
f(\vect{x}) = \bigwedge_{i=0}^{n}\big (\beta_f(S^{\downarrow}_\sigma (i))\wedge x_{\sigma(i)}\big ),
$$
where $x_{\sigma(0)}=0$.
\end{enumerate}
\end{proposition}

\section{Motivating characterizations}\label{sec:MotChar}

Even though horizontal maxitivity and comonotonic maxitivity, as well as their dual counterparts, play an important role in aggregation theory
(as properties shared by noteworthy classes of aggregation functions), they have not yet been described independently. In this section we
investigate each of these properties and determine their corresponding function classes (up to weak versions of monotonicity, in the cases of
horizontal maxitivity and minitivity).

\subsection{Horizontal maxitivity and minitivity}

Recall that a function $f\colon L^{n}\rightarrow L$ is said to be
\begin{itemize}
\item \emph{horizontally maxitive} if, for every $\vect{x}\in L^n$ and every $c\in L$, we have
$$
f(\vect{x}) = f(\vect{x}\wedge c)\vee f([\vect{x}]_c),
$$
where $[\vect{x}]_c$ is the $n$-tuple whose $i$th component is $0$, if $x_i\leqslant c$, and $x_i$, otherwise.

\item \emph{horizontally minitive} if, for every $\vect{x}\in L^n$ and every $c\in L$, we have
$$
f(\vect{x}) = f(\vect{x}\vee c)\wedge f([\vect{x}]^c),
$$
where $[\vect{x}]^c$ is the $n$-tuple whose $i$th component is $1$, if $x_i\geqslant c$, and $x_i$, otherwise.
\end{itemize}
Let us consider the following weak forms of nondecreasing monotonicity:
\begin{itemize}
\item[$\mathbf{(P_1)}$] $f(\vect{e}\wedge c)\leqslant f(\vect{e}'\wedge c)$ for every $\vect{e},\vect{e}'\in\{0,1\}^n$ such that
$\vect{e}\leqslant\vect{e}'$ and every $c\in L$.

\item[$\mathbf{(D_1)}$] $f(\vect{e}\vee c)\leqslant f(\vect{e}'\vee c)$ for every $\vect{e},\vect{e}'\in\{0,1\}^n$ such that
$\vect{e}\leqslant\vect{e}'$ and every $c\in L$.

\item[$\mathbf{(P_2)}$] $f(\vect{e}\wedge c)\leqslant f(\vect{e}\wedge c')$ for every $\vect{e}\in\{0,1\}^n$ and every $c,c'\in L$ such that
$c\leqslant c'$.

\item[$\mathbf{(D_2)}$] $f(\vect{e}\vee c)\leqslant f(\vect{e}\vee c')$ for every $\vect{e}\in\{0,1\}^n$ and every $c,c'\in L$ such that
$c\leqslant c'$.
\end{itemize}

\begin{theorem}\label{thm:HorMax}
A function $f\colon L^n\to L$ is horizontally maxitive and satisfies $\mathbf{P_1}$ if and only if there exists $g\colon L^n\to L$ satisfying
$\mathbf{P_2}$ such that
\begin{equation}\label{eq:HorMaxP1}
f(\vect{x})=\bigvee_{I\subseteq [n]} g\Big(\vect{e}_I\wedge\bigwedge_{i\in I}x_i\Big).
\end{equation}
In this case, we can choose $g=f$.
\end{theorem}

To prove Theorem~\ref{thm:HorMax}, we make use of the following lemma.

\begin{lemma}\label{lemma:DNF-Simplex}
A function $f\colon L^n\to L$ satisfying $\mathbf{P_1}$ is of the form (\ref{eq:HorMaxP1}) if and only if, for every permutation $\sigma$ on
$[n]$ and every $\vect{x}\in L_{\sigma}^n$, we have
$$
f(\vect{x})=\bigvee_{i=1}^{n+1} g\big(\vect{e}_{S_{\sigma}^{\uparrow}(i)}\wedge x_{\sigma(i)}\big).
$$
\end{lemma}

\begin{proof}[Proof of Lemma~\ref{lemma:DNF-Simplex}]
For every permutation $\sigma$ on $[n]$ and every $\vect{x}\in L_{\sigma}^n$, we have
\begin{eqnarray*}
\bigvee_{I\subseteq [n]} g\Big(\vect{e}_I\wedge\bigwedge_{i\in I}x_i\Big) &=& g(\vect{e}_{\varnothing}) \vee \bigvee_{i\in
[n]}\bigvee_{\textstyle{I\subseteq S^{\uparrow}_\sigma (i)\atop \sigma(i)\in I}}g\big(\vect{e}_I\wedge x_{\sigma(i)}\big)\\
&=& g(\vect{e}_{\varnothing}) \vee \bigvee_{i\in [n]}g\big(\vect{e}_{S_{\sigma}^{\uparrow}(i)}\wedge x_{\sigma(i)}\big)\\
&=& \bigvee_{i=1}^{n+1} g\big(\vect{e}_{S_{\sigma}^{\uparrow}(i)}\wedge x_{\sigma(i)}\big).\mbox{\qedhere}
\end{eqnarray*}
\end{proof}

\begin{proof}[Proof of Theorem~\ref{thm:HorMax}]
Let us first show that the condition is sufficient. Let $\sigma$ be a permutation on $[n]$, let $\vect{x}\in L_{\sigma}^n$, $c\in L$, and set
$k=\sup\{i\in [n+1]\colon x_{\sigma(i)}\leqslant c\}$. By Lemma~\ref{lemma:DNF-Simplex}, we have
$$
f(\vect{x}\wedge c) = \bigvee_{i=1}^k g\big(\vect{e}_{S_{\sigma}^{\uparrow}(i)}\wedge x_{\sigma(i)}\big)\vee \bigvee_{i=k+1}^{n+1}
g\big(\vect{e}_{S_{\sigma}^{\uparrow}(i)}\wedge c\big)
$$
and
$$
f([\vect{x}]_c) = \bigvee_{i=k+1}^{n+1} g\big(\vect{e}_{S_{\sigma}^{\uparrow}(i)}\wedge x_{\sigma(i)}\big)
$$
and, since $g$ satisfies $\mathbf{P_2}$, $f$ is horizontal maxitive.

Let us now show that the condition is necessary. Let $\sigma$ be a permutation on $[n]$ and let $\vect{x}\in L_{\sigma}^n$. By repeatedly
applying the horizontal maxitivity with the successive cut levels $x_{\sigma(1)},\ldots,x_{\sigma(n)}$, we obtain
\begin{eqnarray*}
f(\vect{x}) &=& f(\vect{e}_{\{1,\ldots,n\}}\wedge x_{\sigma(1)})\vee f(0,x_{\sigma(2)},\ldots,x_{\sigma(n)})\\
&=& f(\vect{e}_{\{1,\ldots,n\}}\wedge x_{\sigma(1)})\vee f(\vect{e}_{\{2,\ldots,n\}}\wedge x_{\sigma(2)}) \vee
f(0,0,x_{\sigma(3)},\ldots,x_{\sigma(n)})\\
&=& \cdots\\
&=& \bigvee_{i=1}^{n+1} f\big(\vect{e}_{S_{\sigma}^{\uparrow}(i)}\wedge x_{\sigma(i)}\big).
\end{eqnarray*}
Indeed, if for instance $x_{\sigma(1)}=x_{\sigma(2)}<x_{\sigma(3)}$, then
$$
f(x_1,x_2,x_3)=f(\vect{e}_{\{1,2,3\}}\wedge x_{\sigma(1)})\vee f(0,0,x_{\sigma(3)})=f(\vect{e}_{\{1,2,3\}}\wedge x_{\sigma(1)})\vee
f(\vect{e}_{\{3\}}\wedge x_{\sigma(3)})
$$
but, since $f$ satisfies $\mathbf{P_1}$, we have that
$$
f(x_1,x_2,x_3)=f(\vect{e}_{\{1,2,3\}}\wedge x_{\sigma(1)})\vee f(\vect{e}_{\{2,3\}}\wedge x_{\sigma(2)})\vee f(\vect{e}_{\{3\}}\wedge
x_{\sigma(3)}).
$$
Therefore, by Lemma~\ref{lemma:DNF-Simplex}, (\ref{eq:HorMaxP1}) holds with $g=f$. To complete the proof, let us show that $f$ satisfies
$\mathbf{P_2}$. Let $\vect{e}\in\{0,1\}^n$ and let $c,c'\in L$ such that $c\leqslant c'$. Then
$$
f(\vect{e}\wedge c)=f((\vect{e}\wedge c')\wedge c)\leqslant f((\vect{e}\wedge c')\wedge c)\vee f([\vect{e}\wedge c']_c)=f(\vect{e}\wedge
c').\mbox{\qedhere}
$$
\end{proof}

Similarly, we obtain the following dual characterization:

\begin{theorem}\label{thm:HorMin}
A function $f\colon L^n\to L$ is horizontally minitive and satisfies $\mathbf{D_1}$ if and only if there exists $g\colon L^n\to L$ satisfying
$\mathbf{D_2}$ such that
$$
f(\vect{x})=\bigwedge_{I\subseteq [n]} g\Big(\vect{e}_{[n]\setminus I}\vee\bigvee_{i\in I}x_i\Big).
$$
In this case, we can choose $g=f$.
\end{theorem}

From Theorems~\ref{thm:HorMax} and \ref{thm:HorMin} we have the following corollary.

\begin{corollary}\label{cor:HorMaxMin}
A function $f\colon L^n\to L$ is horizontally maxitive (resp.\ horizontally minitive) and satisfies $\mathbf{P_1}$ (resp.\ $\mathbf{D_1}$) if
and only if there are unary nondecreasing functions $\varphi_I\colon L\to L$, for $I\subseteq [n]$, such that
\begin{equation}\label{eq:HorMaxMin}
f(\vect{x})=\bigvee_{I\subseteq [n]}\varphi_I\big(\bigwedge_{i\in I}x_i\big)\qquad \mbox{(resp.\
$\displaystyle{f(\vect{x})=\bigwedge_{I\subseteq [n]}\varphi_I\big(\bigvee_{i\in I}x_i\big)}$)}.
\end{equation}
In this case, we can choose $\varphi_I(x)=f(\vect{e}_I\wedge x)$ (resp.\ $\varphi_I(x)=f(\vect{e}_{[n]\setminus I}\vee x)$) for every
$I\subseteq [n]$.
\end{corollary}

Observe that, by choosing every function $\varphi_I$ in Corollary~\ref{cor:HorMaxMin} as $\varphi_I(x)=f(\vect{e}_I)\wedge f(\vect{e}_I\wedge
x)$ (resp.\ $\varphi_I(x)=f(\vect{e}_{[n]\setminus I})\vee f(\vect{e}_{[n]\setminus I}\vee x)$), equation (\ref{eq:HorMaxMin}) becomes
$$
f(\vect{x})=\bigvee_{I\subseteq [n]}\big(\alpha_f(I)\wedge\bigwedge_{i\in I}\varphi_I(x_i)\big)\qquad \mbox{(resp.\
$\displaystyle{f(\vect{x})=\bigwedge_{I\subseteq [n]}\big(\beta_f(I)\vee\bigvee_{i\in I}\varphi_I(x_i)\big)}$)},
$$
where the set function $\alpha_f$ (resp.\ $\beta_f$) is defined in (\ref{eq:AlphaBeta}).

\begin{remark}
\begin{itemize}
\item[(i)] Theorem~\ref{thm:HorMax} (resp.\ Theorem~\ref{thm:HorMin}) and Corollary~\ref{cor:HorMaxMin} provide descriptions of those
horizontally maxitive (resp.\ horizontally minitive) functions which are nondecreasing.

\item[(ii)] Every Boolean function $f\colon\{0,1\}^n\to\{0,1\}$ satisfying $f(\vect{0})\leqslant f(\vect{x})$ (resp.\ $f(\vect{x})\leqslant
f(\vect{1})$) is horizontally maxitive (resp.\ horizontally minitive). Moreover, not all such functions are nondecreasing, thus showing that
condition $\mathbf{P_1}$ (resp.\ $\mathbf{D_1}$) is necessary in Theorem~\ref{thm:HorMax} (resp.\ Theorem~\ref{thm:HorMin}) and
Corollary~\ref{cor:HorMaxMin}.

\item[(iii)] As shown in \cite{CouMar2}, polynomial functions $f\colon L^n\to L$ are exactly those $\co{\mathcal{R}}_f$-idempotent (i.e.,
satisfying $f(c,\ldots ,c)=c$ for every $c\in \co{\mathcal{R}}_f$) which are nondecreasing, horizontally maxitive, and horizontally minitive.

\item[(iv)] The concept of horizontal maxitivity was introduced, in the case when $L$ is the real interval $[0,1]$, by Benvenuti et
al.~\cite{BenMesViv02} as a general property of the Sugeno integral.
\end{itemize}
\end{remark}

\subsection{Comonotonic maxitivity and minitivity}

Two vectors $\vect{x},\vect{x'}\in L^n$ are said to be \emph{comonotonic} if there exists a permutation $\sigma$ on $[n]$ such that
$\vect{x},\vect{x'}\in L_{\sigma}^n$. A function $f\colon L^{n}\rightarrow L$ is said to be
\begin{itemize}
\item \emph{comonotonic maxitive} if, for any two comonotonic vectors $\vect{x},\vect{x'}\in L^n$, we have
$$
f(\vect{x}\vee \vect{x'}) = f(\vect{x})\vee f(\vect{x'}).
$$

\item \emph{comonotonic minitive} if, for any two comonotonic vectors $\vect{x},\vect{x'}\in L^n$, we have
$$
f(\vect{x}\wedge \vect{x'}) = f(\vect{x})\wedge f(\vect{x'}).
$$
\end{itemize}

Note that for any $\vect{x}\in L^n$ and any $c\in L$, the vectors $\vect{x}\vee c$ and $[\vect{x}]^c$ are comonotonic. As a consequence, if a
function $f\colon L^{n}\rightarrow L$ is comonotonic maxitive (resp.\ comonotonic minitive), then it is horizontally maxitive (resp.\
horizontally minitive). It was also observed in \cite{CouMar2} that if $f$ is comonotonic maxitive or comonotonic minitive, then it is
nondecreasing. Moreover, by using Theorem~\ref{thm:HorMax} and Lemma~\ref{lemma:DNF-Simplex}, we obtain the following result.

\begin{lemma}\label{lemma:ComHor}
A function $f\colon L^{n}\rightarrow L$ is comonotonic maxitive (resp.\ comonotonic minitive) if and only if it is horizontally maxitive (resp.\
horizontally minitive) and satisfies $\mathbf{P_1}$ (resp.\ $\mathbf{D_1}$).
\end{lemma}

\begin{proof}
As observed, the conditions are necessary. We show the sufficiency for comonotonic maxitive functions; the other case follows dually. Let
$f\colon L^{n}\rightarrow L$ be a horizontally maxitive function satisfying $\mathbf{P_1}$ and let $\vect{x},\vect{x}'\in L_{\sigma}^n$ for some
permutation $\sigma$ on $[n]$. By Theorem~\ref{thm:HorMax} and Lemma~\ref{lemma:DNF-Simplex}, there exists $g\colon L^n\to L$ satisfying
$\mathbf{P_2}$ such that
$$
f(\vect{x}\vee\vect{x}') = \bigvee_{i=1}^{n+1} g\big(\vect{e}_{S_{\sigma}^{\uparrow}(i)}\wedge (x_{\sigma(i)}\vee x'_{\sigma(i)})\big).
$$
By distributivity and $\mathbf{P_2}$, we have
\begin{eqnarray*}
f(\vect{x}\vee\vect{x}') &=& \bigvee_{i=1}^{n+1} g\big((\vect{e}_{S_{\sigma}^{\uparrow}(i)}\wedge x_{\sigma(i)})\vee
(\vect{e}_{S_{\sigma}^{\uparrow}(i)}\wedge
x'_{\sigma(i)})\big)\\
&=& \bigvee_{i=1}^{n+1} \Big(g\big(\vect{e}_{S_{\sigma}^{\uparrow}(i)}\wedge x_{\sigma(i)}\big)\vee
g\big(\vect{e}_{S_{\sigma}^{\uparrow}(i)}\wedge
x'_{\sigma(i)}\big)\Big)%
~=~ f(\vect{x})\vee f(\vect{x}').\mbox{\qedhere}
\end{eqnarray*}
\end{proof}

Combining Theorems~\ref{thm:HorMax} and \ref{thm:HorMin} with Lemma~\ref{lemma:ComHor}, we immediately obtain the descriptions of the classes of
comonotonic maxitive and comonotonic minitive functions.

\begin{theorem}\label{thm:ComMax}
A function $f\colon L^n\to L$ is comonotonic maxitive if and only if there exists $g\colon L^n\to L$ satisfying $\mathbf{P_2}$ such that
$$
f(\vect{x})=\bigvee_{I\subseteq [n]} g\Big(\vect{e}_I\wedge\bigwedge_{i\in I}x_i\Big).
$$
In this case, we can choose $g=f$.
\end{theorem}

\begin{theorem}\label{thm:ComMin}
A function $f\colon L^n\to L$ is comonotonic minitive if and only if there exists $g\colon L^n\to L$ satisfying $\mathbf{D_2}$ such that
$$
f(\vect{x})=\bigwedge_{I\subseteq [n]} g\Big(\vect{e}_{[n]\setminus I}\vee\bigvee_{i\in I}x_i\Big).
$$
In this case, we can choose $g=f$.
\end{theorem}

As before, we have the following corollary.

\begin{corollary}\label{cor:ComMaxMin}
A function $f\colon L^n\to L$ is comonotonic maxitive (resp.\ comonotonic minitive) if and only if there are unary nondecreasing functions
$\varphi_I\colon L\to L$, for $I\subseteq [n]$, such that
$$
f(\vect{x})=\bigvee_{I\subseteq [n]}\varphi_I\big(\bigwedge_{i\in I}x_i\big)\qquad \mbox{(resp.\
$\displaystyle{f(\vect{x})=\bigwedge_{I\subseteq [n]}\varphi_I\big(\bigvee_{i\in I}x_i\big)}$)}.
$$
In this case, we can choose $\varphi_I(x)=f(\vect{e}_I\wedge x)$ (resp.\ $\varphi_I(x)=f(\vect{e}_{[n]\setminus I}\vee x)$) for every
$I\subseteq [n]$.
\end{corollary}

\begin{remark}
\begin{itemize}
\item[(i)] An alternative description of comonotonic maxitive (resp.\ comonotonic minitive) functions was obtained in Grabisch et
al.~\cite[Chapter 2]{GraMarMesPap09} in the case when $L$ is a real interval.

\item[(ii)] It was shown in \cite{CouMar2} that polynomial functions $f\colon L^n\to L$ are exactly those $\co{\mathcal{R}}_f$-idempotent
functions which are comonotonic maxitive and comonotonic minitive.

\item[(ii)] Comonotonic minitivity and maxitivity were introduced in the context of Sugeno integrals in de Campos et al.~\cite{deCLamMor91}.
\end{itemize}
\end{remark}

\section{Quasi-polynomial functions}

Motivated by the results of Section~\ref{sec:MotChar} concerning horizontal maxitivity and comonotonic maxitivity, as well as their dual
counterparts, we now study combinations of these properties. This will lead to a relaxation of the notion of polynomial function, which we will
refer to as \emph{quasi-polynomial function}. Accordingly, we introduce weaker variants of well-established properties, such as homogeneity and
median decomposability, which are then used to provide further axiomatizations of the class of quasi-polynomial functions.

\subsection{Motivation and definition}

We start by looking at combinations of those properties studied in Section~\ref{sec:MotChar}. These are considered in the following result.

\begin{theorem}\label{thm:QuasiPol}
Let $f\colon L^n\to L$ be a function. The following assertions are equivalent:
\begin{itemize}
\item[(i)] $f$ is horizontally maxitive, horizontally minitive, and satisfies $\mathbf{P_1}$ or $\mathbf{D_1}$.

\item[(ii)] $f$ is comonotonic maxitive and comonotonic minitive.

\item[(iii)] $f$ is horizontally maxitive and comonotonic minitive.

\item[(iv)] $f$ is comonotonic maxitive and horizontally minitive.

\item[(v)] There exist a polynomial function $p\colon L^n\to L$ and a nondecreasing function $\varphi\colon L\to L$ such that
$$
f(x_1,\ldots,x_n)=p(\varphi(x_1),\ldots,\varphi(x_n)).
$$
\end{itemize}
If these conditions hold then we can choose for $p$ the unique polynomial function $p_f$ extending $f|_{\{0,1\}^n}$ and for $\varphi$ the
diagonal section $\delta_f$ of $f$.
\end{theorem}

\begin{proof}
The equivalences between assertions $(i)$--$(iv)$ follow from Lemma~\ref{lemma:ComHor}. Note that only one of the conditions $\mathbf{P_1}$ and
$\mathbf{D_1}$ suffices in assertion $(i)$ since the other one then follows from Theorem~\ref{thm:HorMax} or Theorem~\ref{thm:HorMin}, which
ensure nondecreasing monotonicity. Also, $\mathbf{P_1}$ (resp.\ $\mathbf{D_1}$) is not needed in assertion $(iii)$ (resp.\ $(iv)$) since, as
already observed, comonotonic minitivity (resp.\ comonotonic maxitivity) ensures nondecreasing monotonicity. To see that $(v)\Rightarrow (ii)$
holds, just note that every polynomial function is comonotonic maxitive and comonotonic minitive, and that if $\vect{x},\vect{x'}\in
L_{\sigma}^n$ for some permutation $\sigma$ on $[n]$ then $\varphi(\vect{x}),\varphi(\vect{x'})\in L_{\sigma}^n$. To conclude the proof of the
theorem, it is enough to show that $(iii)\Rightarrow (v)$. By Theorem~\ref{thm:HorMax} and Lemma~\ref{lemma:DNF-Simplex}, for every permutation
$\sigma$ on $[n]$ and every $\vect{x}\in L_{\sigma}^n$, we have
$$
f(\vect{x}) = \bigvee_{i=1}^{n+1} f\big(\vect{e}_{S_{\sigma}^{\uparrow}(i)}\wedge x_{\sigma(i)}\big).
$$
Since the vectors $\vect{e}_{S_{\sigma}^{\uparrow}(i)}$ and $(x_{\sigma(i)},\ldots,x_{\sigma(i)})$ are comonotonic, we get
\begin{eqnarray*}
f(\vect{x}) &=& \bigvee_{i=1}^{n+1} \Big(f\big(\vect{e}_{S_{\sigma}^{\uparrow}(i)}\big)\wedge \delta_f(x_{\sigma(i)})\Big) =
p(\varphi(x_1),\ldots,\varphi(x_n)),
\end{eqnarray*}
where $p$ is the unique polynomial function $p_f$ extending $f|_{\{0,1\}^n}$ (which exists due to $\mathbf{P_1}$) and $\varphi$ is the diagonal
section $\delta_f$ of $f$.
\end{proof}

Theorem~\ref{thm:QuasiPol} motivates the following definition.

\begin{definition}
We say that a function $f\colon L^n\to L$ is a \emph{quasi-polynomial function} (resp.\ a \emph{discrete quasi-Sugeno integral}, a
\emph{quasi-term function}) if there exist a polynomial function (resp.\ a discrete Sugeno integral, a term function) $p\colon L^n\to L$ and a
nondecreasing function $\varphi\colon L\to L$ such that $f=p\circ\varphi$, that is,
\begin{equation}\label{eq:QuasiPol2}
f(x_1,\ldots,x_n)=p(\varphi(x_1),\ldots,\varphi(x_n)).
\end{equation}
\end{definition}

\begin{remark}\label{remark:ND}
\begin{itemize}
\item[(i)] Note that each quasi-polynomial function $f\colon L^n\to L$ can be represented as a combination of constants and a nondecreasing
unary function $\varphi$ (applied to the projections $\vect{x}\mapsto x_i$) using the lattice operations $\vee $ and $\wedge$.

\item[(ii)] In the setting of decision-making under uncertainty, the nondecreasing function $\varphi$ in (\ref{eq:QuasiPol2}) can be thought of
as a \emph{utility function} and the corresponding quasi-polynomial function as a (qualitative) \emph{global preference functional}; see for
instance Dubois et al.~\cite{DubMarPraRouSab01}.
\end{itemize}
\end{remark}

Note that the functions $p$ and $\varphi$ in (\ref{eq:QuasiPol2}) are not necessarily unique. For instance, if $f$ is a constant $c\in L$, then
we could choose $p\equiv c$ and $\varphi$ arbitrarily, or $p$ idempotent and $\varphi\equiv c$. To describe all possible choices for $p$ and
$\varphi$, we shall make use of the following result implicit in \cite{CouMar2}. For any integers $m,n\geqslant 1$, any vector $\vect{x}\in
L^m$, and any function $f\colon L^{n}\rightarrow L$, we define $\langle\vect{x}\rangle_f\in L^m$ as the $m$-tuple
$$
\langle\vect{x}\rangle_f=\median(f(\vect{0}),\vect{x},f(\vect{1})),
$$
where the right-hand side median is taken componentwise.

\begin{lemma}\label{lemma:45474786}
Every polynomial function $p\colon L^n\to L$ satisfies
$$
p(\vect{x}\vee c)=p(\vect{x})\vee\langle c\rangle_p \quad\mbox{and}\quad p(\vect{x}\wedge c)=p(\vect{x})\wedge\langle c\rangle_p
$$
for every $\vect{x}\in L^n$ and every $c\in L$.
\end{lemma}

\begin{proposition}\label{prop:DescrPP564}
Let $f\colon L^n\to L$ be a quasi-polynomial function and let $p_f:L^n\to L$ be the unique polynomial function extending $f|_{\{0,1\}^n}$. We
have
$$
\{(p,\varphi)\colon f=p\circ\varphi\} = \{(p,\varphi)\colon p_f=\langle p\rangle_f~\mbox{and}~\delta_f=\langle\varphi\rangle_p\},
$$
where $p$ and $\varphi$ stand for polynomial and unary nondecreasing functions, respectively. In particular, we have $f=p_f\circ\delta_f$.
\end{proposition}

\begin{proof}
$(\subseteq)$ Let $p$ and $\varphi$ be such that $f=p\circ\varphi$. First observe that, for any $c\in L$, we have
$\delta_f(c)=(\delta_p\circ\varphi)(c)=\langle\varphi(c)\rangle_p$. By assertion (b) of Proposition~\ref{Goodstein}, to complete the proof, it
is enough to show the equality $p_f=\langle p\rangle_f$ restricted to $\{0,1\}^n$.

By Lemma~\ref{lemma:45474786}, for any $\vect{e}\in\{0,1\}^n$we have
\begin{eqnarray*}
p_f(\vect{e}) &=& p\big((\varphi(0)\vee\vect{e})\wedge\varphi(1)\big) ~=~ \big(\langle\varphi(0)\rangle_p\vee p(\vect{e})\big)\wedge
\langle\varphi(1)\rangle_p\\
&=& \big(\delta_f(0)\vee p(\vect{e})\big)\wedge \delta_f(1) ~=~ \langle p(\vect{e})\rangle_f.
\end{eqnarray*}

$(\supseteq)$ Let $p$ and $\varphi$ be such that $p_f=\langle p\rangle_f$ and $\delta_f=\langle\varphi\rangle_p$. Again we have
\begin{eqnarray*}
p_f(\vect{e}) &=& \langle p(\vect{e})\rangle_f ~=~ \big(\delta_f(0)\vee p(\vect{e})\big)\wedge \delta_f(1)
~=~ \big(\langle\varphi(0)\rangle_p\vee p(\vect{e})\big)\wedge \langle\varphi(1)\rangle_p\\
&=& \big\langle p\big((\varphi(0)\vee\vect{e})\wedge\varphi(1)\big)\big\rangle_p ~=~ (p\circ\varphi)(\vect{e}).\qedhere
\end{eqnarray*}
\end{proof}

It was shown in \cite{Marc} that every polynomial function $p\colon L^n\to L$ can be represented as $\langle q\rangle_p$ for some
discrete Sugeno integral $q\colon L^n\to L$. Combining this with Proposition~\ref{prop:DescrPP564}, we obtain the next result.

\begin{corollary}\label{cor:Sug53453}
The class of quasi-polynomial functions is exactly the class of discrete quasi-Sugeno integrals.
\end{corollary}

\subsection{Further axiomatizations}

We now recall some properties of polynomial functions, namely homogeneity and median decomposability, and we propose weaker variants of these to
provide alternative axiomatizations of the class of quasi-polynomial functions.

\subsubsection{Quasi-homogeneity}

Let $S$ be a subset of $L$. A function $f\colon L^{n}\rightarrow L$ is said to be \emph{$S$-max homogeneous} (resp.\ \emph{$S$-min homogeneous})
if for every $\vect{x}\in L^n$ and every $c\in S$, we have
$$
f(\vect{x}\vee c) = f(\vect{x})\vee c \qquad \textrm{(resp.\ $f(\vect{x}\wedge c) = f(\vect{x})\wedge c$}).
$$
Although polynomial functions $p\colon L^n\to L$ share both of these properties for any $S\subseteq\co{\mathcal{R}}_p$, this is not the case for
quasi-polynomial functions. For instance, let $f_1,f_2\colon [0,1]\rightarrow [0,1]$  be respectively given by $f_1(x)=x^2$ and
$f_2(x)=\sqrt{x}$. Clearly, $f_1$ and $f_2$ are quasi-polynomial functions but, e.g., for $x=c\in\left]0,1\right[$, we have
$$
f_1(x\vee c) < f_1(x)\vee c\quad\mbox{and}\quad f_2(x\wedge c) > f_2(x)\wedge c.
$$
This example motivates the following relaxations. We say that a function $f\colon L^{n}\rightarrow L$ is \emph{quasi-max homogeneous} (resp.\
\emph{quasi-min homogeneous}) if for every $\vect{x}\in L^n$ and $c\in L$, we have
$$
f(\vect{x}\vee c) = f(\vect{x})\vee  \delta_f(c) \qquad \textrm{(resp.\ $f(\vect{x}\wedge c) = f(\vect{x})\wedge \delta_f(c)$}).
$$
Observe that if $f$ is $\co{\mathcal{R}}_f$-idempotent (i.e., satisfying $f(c,\ldots ,c)=c$ for every $c\in \co{\mathcal{R}}_f$), then
$\co{\mathcal{R}}_f$-min homogeneity (resp.\ $\co{\mathcal{R}}_f$-max homogeneity) is equivalent to quasi-min homogeneity (resp.\ quasi-max
homogeneity).


\begin{lemma}\label{quasi-Hor-Min-Hom}
Let $f\colon L^{n}\rightarrow L$ be nondecreasing and quasi-min homogeneous (resp.\ quasi-max homogeneous). Then $f$ is quasi-max homogeneous
(resp.\ quasi-min homogeneous) if and only if it is horizontally maxitive (resp.\ horizontally minitive).
\end{lemma}

\begin{proof}
Let $f\colon L^{n}\rightarrow L$ be nondecreasing and quasi-min homogeneous and suppose first that
 $f$ is also quasi-max homogeneous. For any $\vect{x}\in L^n$ and any $c\in L$,
we have
\begin{eqnarray*}
  f(\vect{x}\wedge c)\vee f([\vect{x}]_c)
  &=& \big(f(\vect{x})\wedge \delta_f(c)\big)\vee f([\vect{x}]_c) ~=~ \big(f(\vect{x})\vee f([\vect{x}]_c)\big)\wedge\big(\delta_f(c)\vee f([\vect{x}]_c)\big) \\
  &=& f(\vect{x})\wedge f(c\vee [\vect{x}]_c) ~=~ f(\vect{x}).
\end{eqnarray*}
Therefore, $f$ is horizontally maxitive.

Now assume that $f$ is horizontally maxitive. From horizontal maxitivity and quasi-min homogeneity, it follows that for any $\vect{x}\in L^n$
and any $c\in L$, we have $ f(\vect{x}\vee c) = \delta_f(c) \vee f([\vect{x}]_c). $ Thus,
$$
f(\vect{x}\vee c) = f(\vect{x}\wedge c)\vee f(\vect{x}\vee c) = f(\vect{x}\wedge c) \vee f([\vect{x}]_c) \vee \delta_f(c) = f(\vect{x})\vee
\delta_f(c)
$$
and hence $f$ is quasi-max homogeneous. The remaining claim can be verified dually.
\end{proof}

Combining Theorem~\ref{thm:QuasiPol} and Lemma~\ref{quasi-Hor-Min-Hom}, we obtain a characterization of quasi-polynomial functions in terms of
quasi-min homogeneity and quasi-max homogeneity. However, we provide a constructive proof.

\begin{theorem}\label{theorem:Quasi-Hom}
 A function $f\colon L^{n}\rightarrow L$ is a quasi-polynomial function if and only if it is nondecreasing, quasi-max homogeneous, and quasi-min homogeneous.
\end{theorem}

\begin{proof}
The necessity of the conditions follows from Proposition~\ref{prop:DescrPP564} together with the facts that
$\mathcal{R}_{\delta_f}\subseteq\mathcal{R}_f\subseteq \mathcal{R}_{p_f}$ and that $p_f$ is both $\mathcal{R}_{p_f}$-max homogeneous and
$\mathcal{R}_{p_f}$-min homogeneous.

To verify the sufficiency, let $\vect{x}\in L^n$. By nondecreasing monotonicity and quasi-min homogeneity, for every $I\subseteq [n]$ we have
$$
f(\vect{x}) \geqslant f\big(\vect{e}_I \wedge \bigwedge_{i\in I}x_i\big) = f(\vect{e}_I) \wedge \delta_f\big(\bigwedge_{i\in I}x_i\big)
$$
and thus $f(\vect{x})\geqslant \bigvee_{I\in [n]}  f(\vect{e}_I) \wedge\delta_f(\bigwedge_{i\in I}x_i)$. To complete the proof, it is enough to
establish the converse inequality since $\delta_f(\bigwedge_{i\in I}x_i)=\bigwedge_{i\in I}\delta_f (x_i)$. Let $I^*\subseteq [n]$ be such that
$f(\vect{e}_{I^*}) \wedge\delta_f( \bigwedge_{i\in I^*}x_i)$ is maximum. Define
$$
J=\big\{j\in [n]\colon \delta_f(x_j)\leqslant f(\vect{e}_{I^*}) \wedge\delta_f\big( \bigwedge_{i\in I^*}x_i\big)\big\}.
$$
We claim that $J\neq\varnothing$. For the sake of contradiction, suppose that $\delta_f(x_j)> f(\vect{e}_{I^*}) \wedge\delta_f(\bigwedge_{i\in
I^*}x_i)$ for every $j\in [n]$. Then, by nondecreasing monotonicity, we have $f(\vect{e}_{[n]})\geqslant f(\vect{e}_{I^*})$, and since
$f(\vect{e}_{[n]})= \delta_f(1)\geqslant\delta_f(\bigwedge_{i\in [n]}x_i)$,
$$
f(\vect{e}_{[n]}) \wedge \delta_f\big(\bigwedge_{i\in [n]}x_i\big) > f(\vect{e}_{I^*}) \wedge\delta_f\big(\bigwedge_{i\in I^*}x_i\big),$$ which
contradicts the definition of $I^*$. Thus $J\neq \varnothing$.

Now, let $j\in J$ such that $x_j=\sup \{x_k\colon k\in J\}$. By nondecreasing monotonicity and quasi-max homogeneity, we have
$$
f(\vect{x}) \leqslant f(x_j\vee \vect{e}_{[n]\setminus J}) = \delta_f(x_j)\vee f(\vect{e}_{[n]\setminus J})\leqslant \big( f(\vect{e}_{I^*})
\wedge\delta_f\big( \bigwedge_{i\in I^*}x_i\big)\big)\vee f(\vect{e}_{[n]\setminus J}).
$$
We claim that $f(\vect{e}_{[n]\setminus J})\leqslant f(\vect{e}_{I^*}) \wedge \delta_f( \bigwedge_{i\in I^*}x_i)$. Otherwise, by definition of
$J$ we would have
$$
f(\vect{e}_{[n]\setminus J})\wedge \delta_f\big(\bigwedge_{i\in [n]\setminus J}x_i\big) > f(\vect{e}_{I^*}) \wedge \delta_f\big( \bigwedge_{i\in
I^*}x_i\big),
$$
which would again contradict the choice of $I^*$. Thus,
$$
f(\vect{x}) \leqslant f(\vect{e}_{I^*}) \wedge \delta_f\big( \bigwedge_{i\in I^*}x_i\big) = \bigvee_{I\in [n]} \big(f(\vect{e}_I) \wedge
\delta_f\big(  \bigwedge_{i\in I}x_i\big)\big).\qedhere
$$
\end{proof}

\subsubsection{Quasi-median decomposability}

A function $f\colon L^n\to L$ is said to be \emph{median decomposable} if, for every $\vect{x}\in L^n$ and every $k\in [n]$, we have
$$
f(\vect{x})=\median\big(f(\vect{x}^{0}_{k}), x_k, f(\vect{x}^{1}_{k})\big),
$$
where $\vect{x}^{c}_{k} = (x_1,\ldots, x_{k-1},c,x_{k+1},\ldots ,x_n)$ for any $c\in L$. As shown in \cite{Marc}, the class of polynomial
functions are exactly those functions which are median decomposable.

In complete analogy with the previous subsection we propose the following weaker variant of median decomposability. We say that a function
$f\colon L^{n}\rightarrow L$ is \emph{quasi-median decomposable} if, for every $\vect{x}\in L^n$ and every $k\in [n]$, we have
$$
f(\vect{x})=\median\big(f(\vect{x}^{0}_{k}), \delta_f(x_k), f(\vect{x}^{1}_{k})\big).
$$
Note that every nondecreasing unary function is quasi-median decomposable.
%


Observe that $\vee$ and $\wedge$, as well as any nondecreasing function $\varphi\colon L\to L$, are quasi-median decomposable. Also, it is easy
to see that any combination of constants and a nondecreasing unary function $\varphi$ using $\vee$ and $\wedge$ is quasi-median decomposable and
hence, by Remark~\ref{remark:ND} (i), every quasi-polynomial function is quasi-median decomposable. Our following result asserts that
quasi-median decomposable functions $f\colon L^{n}\rightarrow L$ with a nondecreasing diagonal section $\delta_f$ are exactly the
quasi-polynomial functions.

\begin{theorem}\label{theorem:QuasiMedian}
A function $f\colon L^{n}\rightarrow L$ is a quasi-polynomial function if and only if
 $\delta_f$ is nondecreasing and $f$ is quasi-median decomposable.
\end{theorem}

\begin{proof}
By the above observation, we only need to verify that the conditions are sufficient. By Theorem~\ref{theorem:Quasi-Hom}, it is enough to show
that if
 $\delta_f$ is nondecreasing and $f$ is quasi-median decomposable, then $f$ is nondecreasing, quasi-max homogeneous, and quasi-min homogeneous.
Since nondecreasing monotonicity can be verified on vectors differing only on a single component, by assuming that $\delta_f$ is nondecreasing
and $f$ is quasi-median decomposable, it follows that $f$ is nondecreasing.

We show that $f$ is quasi-min homogeneous. The dual property follows similarly. Let $\vect{x}\in L^n$ and $c\in L$. Clearly, we have
$f(\vect{x}\wedge c)=f(\vect{x})\wedge \delta_f(c)$ whenever $c<\bigwedge_{i\in[n]}x_i$ or $c>\bigvee_{i\in[n]}x_i$. So suppose that
$\bigwedge_{i\in[n]}x_i\leqslant c\leqslant \bigvee_{i\in[n]}x_i$ and, without loss of generality, assume that $x_1\leqslant \cdots \leqslant
x_n$.

\begin{claim}\label{claim1}
If $c\leqslant x_j$ then $f(x_1,\ldots ,x_j, 1, \ldots ,1)\wedge \delta_f(c)= f(x_1,\ldots ,x_{j-1}, 1, \ldots ,1)\wedge \delta_f(c)$.
\end{claim}

\begin{proof}[Proof of Claim~\ref{claim1}]
By nondecreasing monotonicity and quasi-median decomposability,
$$
f(x_1,\ldots ,x_j, 1, \ldots ,1)=\big( f(x_1,\ldots ,x_{j-1}, 0,1, \ldots ,1)\vee \delta_f(x_j)\big) \wedge f(x_1,\ldots ,x_{j-1}, 1, \ldots
,1).
$$
Since $\delta_f(c)\leqslant \delta_f(x_j)$, the claim follows.
\end{proof}

Let $k=\sup\{i\in [n]\colon x_i<c\}$. By repeated applications of Claim~\ref{claim1}, it follows that
\begin{eqnarray*}
f(\vect{x})\wedge \delta_f(c) &=& f(x_1,\ldots ,x_{k}, 1, \ldots ,1)\wedge \delta_f(c)\\
&=& f(x_1,\ldots ,x_{k}, c, \ldots ,c)\wedge \delta_f(c) ~=~ f(\vect{x}\wedge c)\wedge \delta_f(c).
\end{eqnarray*}
Moreover, by nondecreasing monotonicity, we have $f(\vect{x}\wedge c)=f(\vect{x}\wedge c)\wedge \delta_f(c)$, and thus $f(\vect{x}\wedge
c)=f(\vect{x})\wedge \delta_f(c)$.
\end{proof}

\section{Some special classes of quasi-polynomial functions}

In this final section we consider few noteworthy subclasses of quasi-polynomial functions, namely those of quasi-term functions and
quasi-weighted maximum and minimum functions, and provide characterizations accordingly.

\subsection{Quasi-term functions}

We say that a function $f\colon L^n\to L$ is
\begin{itemize}
\item \emph{conservative} if, for every $\vect{x}\in L^n$, we have $f(\vect{x})\in \{x_1,\ldots,x_n\}$.

\item \emph{quasi-conservative} if, for every $\vect{x}\in L^n$, we have $f(\vect{x})\in \{\delta_f(x_1),\ldots,\delta_f(x_n)\}$.
\end{itemize}
Note that, if $f$ is idempotent, then it is quasi-conservative if and only if it is conservative.

\begin{theorem}
A quasi-polynomial function $f\colon L^n\to L$ is a quasi-term function if and only if it is quasi-conservative.
\end{theorem}

\begin{proof}
The condition is clearly necessary. To show that it is also sufficient, we use Proposition~\ref{prop:DescrPP564} and note that
$$
f(\vect{x})=(p_f\circ\delta_f)(\vect{x})=\bigvee_{\textstyle{I\subseteq [n]\atop f(\vect{e}_I)=1}} \bigwedge_{i\in
I}\delta_f(x_i)=\delta_f\Bigg(\bigvee_{\textstyle{I\subseteq [n]\atop f(\vect{e}_I)=1}} \bigwedge_{i\in I}x_i\Bigg).\mbox{\qedhere}
$$
\end{proof}

\subsection{Quasi-weighted maximum and minimum functions}

A function $f\colon L^{n}\rightarrow L$ is said to be a \emph{weighted maximum function} if there are $v_0,v_1,\ldots ,v_n\in L$ such that
\begin{equation}\label{eq:Wmax}
f(\vect{x})= v_0\vee \bigvee_{i\in [n]} (v_i\wedge x_i).
\end{equation}
Similarly, $f\colon L^{n}\rightarrow L$ is said to be a \emph{weighted minimum function} if there are $w_0,w_1,\ldots ,w_n\in L$ such that
\begin{equation}\label{eq:Wmin}
f(\vect{x})= w_0\wedge \bigwedge_{i\in [n]} (w_i\vee x_i).
\end{equation}
We say that a function $f\colon L^n\to L$ is a \emph{quasi-weighted maximum function} (resp.\ a \emph{quasi-weighted minimum function}) if there
exist a weighted maximum function (resp.\ a weighted minimum function) $p\colon L^n\to L$ and a nondecreasing function $\varphi\colon L\to L$
such that $f=p\circ\varphi$.

To present a axiomatization of each of these classes, we need to recall some terminology. We say that a function $f\colon L^n\to L$ is
\begin{itemize}
\item \emph{maxitive} if, for every $\vect{x},\vect{x}'\in L^n$, we have $f(\vect{x}\vee\vect{x}')=f(\vect{x})\vee f(\vect{x}')$.

\item \emph{minitive} if, for every $\vect{x},\vect{x}'\in L^n$, we have $f(\vect{x}\wedge\vect{x}')=f(\vect{x})\wedge f(\vect{x}')$.
\end{itemize}

We first recall the descriptions of maxitive and minitive functions; see \cite{DubPra90,Mar00c}. For the
sake of self-containment, a short proof is given here.

\begin{proposition}\label{prop:13137}
A function $f\colon L^n\to L$ is maxitive (resp.\ minitive) if and only if there are nondecreasing unary functions $f_i\colon L\to L$ $(i\in
[n])$ such that, for every $\vect{x}\in L^n$,
$$
f(\vect{x})=\bigvee_{i\in [n]} f_i(x_i)\qquad (\mbox{resp.}~f(\vect{x})=\bigwedge_{i\in [n]} f_i(x_i)).
$$
\end{proposition}

\begin{proof}
Sufficiency follows from the fact that the functions $f_i$ are nondecreasing. Let us show the necessity for maxitive functions only; the other
case follows dually. By maxitivity, we have
\begin{equation}\label{eq:fVeeMax}
f(\vect{x})=\bigvee_{i\in [n]} f(\vect{0}_i^{x_i})=\bigvee_{i\in [n]} f_i(x_i),
\end{equation}
where, for every $i\in [n]$, $f_i\colon L\to L$ is defined by $f_i(x)=f(\vect{0}_i^{x})$. To see that each $f_i$ is nondecreasing, just observe
that each is maxitive.
\end{proof}

\begin{theorem}
Let $f\colon L^n\to L$ be a quasi-polynomial function. Then $f$ is a quasi-weighted maximum function (resp.\ quasi-weighted minimum function) if
and only if it is maxitive (resp.\ minitive).
\end{theorem}

\begin{proof}
By Proposition~\ref{prop:13137}, the condition is clearly necessary. We show the sufficiency for maxitive functions; the other case follows
dually. Assume $f$ is a maxitive quasi-polynomial function. By definition, there exist a polynomial function $p\colon L^n\to L$ and a
nondecreasing function $\varphi\colon L\to L$ such that $f=p\circ\varphi$. Then, by Lemma~\ref{lemma:45474786} and
Proposition~\ref{prop:DescrPP564}, for every $i\in [n]$, we have
\begin{eqnarray*}
f(\vect{0}_i^{x_i}) &=& p\big(\varphi(0)\vee\big(\vect{e}_{\{i\}}\wedge\varphi(x_i)\big)\big)%
~=~ \langle\varphi(0)\rangle_p \vee\big(p(\vect{e}_{\{i\}})\wedge\langle\varphi(x_i)\rangle_p\big)\\
&=& \delta_f(0) \vee\big(p(\vect{e}_{\{i\}})\wedge\delta_f(x_i)\big).
\end{eqnarray*}
Setting $v_0=\delta_f(0)$ and $v_i=p(\vect{e}_{\{i\}})$ for $i=1,\ldots,n$, by (\ref{eq:fVeeMax}), we finally obtain
$$
f(\vect{x})=v_0\vee\bigvee_{i\in [n]}\big(v_i\wedge\delta_f(x_i)\big).\qedhere
$$
\end{proof}

\begin{remark}
\begin{itemize}
\item[(i)] Idempotent weighted maximum functions $f\colon L^{n}\rightarrow L$ are those functions (\ref{eq:Wmax}) for which $v_0=0$ and
$\vee_{i\in [n]} v_i=1$. Dually, idempotent weighted minimum functions $f\colon L^{n}\rightarrow L$ are those functions (\ref{eq:Wmin}) for
which $w_0=1$ and $\wedge_{i\in [n]} w_i=0$. These functions were introduced on real intervals by Dubois and Prade~\cite{DubPra86} in fuzzy set
theory.

\item[(ii)] As observed in Proposition~\ref{prop:DescrPP564}, the underlying weighted maximum function (resp.\ weighted minimum function)
defining a given quasi-weighted maximum function (resp.\ quasi-weighted minimum function) can be chosen to be idempotent.
\end{itemize}
\end{remark}


\end{document}